\documentclass[12pt,reqno]{amsart}
\usepackage{amsmath, amssymb, amsfonts, xcolor}
\usepackage{hyperref}
\numberwithin{equation}{section}
\usepackage{amsthm}
\newtheorem{thm}{Theorem}[section]

\newtheorem{prop}[thm]{Proposition}
\newtheorem{lem}[thm]{Lemma}
\newtheorem{conj}[thm]{Conjecture}
\newtheorem{dfn}[thm]{Definition}

\newtheorem{remark}[thm]{Remark}
\newtheorem{cor}[thm]{Corollary}

\usepackage[a4paper, top = 1.2in, bottom = 1.2in, left = 1.2in, right = 1.2in]{geometry}

\usepackage{enumitem}
\newlist{steps}{enumerate}{1}
\setlist[steps, 1]{label = Step \arabic*:}

\numberwithin{equation}{section}

\newcommand{\F}{\mathbb{F}}
\newcommand{\N}{\mathbb{N}}
\newcommand{\Q}{\mathbb{Q}}

\newcommand{\Z}{\mathbb{Z}}

\newcommand{\f}{\mathbf{f}}

\newcommand{\mcO}{\mathcal{O}}

\newcommand{\mfa}{\mathfrak{a}}

\newcommand{\mfc}{\mathfrak{c}}

\newcommand{\mfm}{\mathfrak{m}}
\newcommand{\mfn}{\mathfrak{n}}

\newcommand{\mfp}{\mathfrak{p}}
\newcommand{\mfq}{\mathfrak{q}}
\newcommand{\mfP}{\mathfrak{P}}

\newcommand{\GL}{\mathrm{GL}}

\newcommand{\Gal}{\mathrm{Gal}}

\def\1{1\!\!1}
\newcommand{\mrm}[1]{\mathrm{#1}}

\title[On the solutions of generalized Fermat equation over $K$]{On the solutions of the generalized Fermat equation over totally real number fields}

\author[S. Sahoo]{Satyabrat Sahoo}
\address[S. Sahoo]{Yau Mathematical Sciences Center, Tsinghua University,\\
	Beijing 100084, China.}
\email{sahoos@mail.tsinghua.edu.cn}

\date{\today}

\begin{document}
	\begin{abstract}
		Let $K$ be a totally real number field and $ \mcO_K$ be the ring of integers of $K$. In this article, we study the asymptotic solutions of the generalized Fermat equation, namely $Ax^p+By^p+Cz^p=0$ over $K$ with prime exponent $p$, where $A,B,C \in \mcO_K \setminus \{0\}$.
		For a certain class of fields $K$, we prove that the equation $Ax^p+By^p+Cz^p=0$ has no asymptotic solution $(a,b,c) \in \mcO_K^3$ with $2\mid abc$. Then, under some assumptions on $A,B,C$, we also prove that $Ax^p+By^p+Cz^p=0$ has no asymptotic solution in $K^3$.
		%
		Finally, we give several purely local criteria for the nonexistence of asymptotic solutions in $K^3$. 
	\end{abstract}
	\maketitle
	\section{Introduction}
	The study of Diophantine equations is a longstanding topic in number theory. The well-known Fermat equation is an excellent example. Wiles~\cite{W95} proved that the Fermat equation $x^n+y^n=z^n$ has no non-trivial integer solutions for integers $n \geq 3$. The proof relies on the modularity of semistable elliptic curves $E$ over $\Q$ (cf. \cite{W95,TW95}), the irreducibility of the mod-$p$ Galois representations $\bar{\rho}_{E,p}$ attached to $E$ (cf. \cite{M78}), and Ribet's level-lowering theorem for $\bar{\rho}_{E,p}$ (cf. \cite{R90}). 
	
	After Wiles proved Fermat's Last Theorem, there was significant progress in the study of the generalized Fermat equation, namely
	\begin{equation}
		\label{generalized Fermat eqn}
		Ax^p+By^q+Cz^r=0,
	\end{equation}
	where $A,B,C$ are nonzero coprime integers and $p,q,r \geq 2$ satisfy $\frac{1}{p}+\frac{1}{q}+\frac{1}{r}<1$. We call $(p,q,r)$ the signature of~\eqref{generalized Fermat eqn}.
	The following conjecture is known for the generalized Fermat equation~\eqref{generalized Fermat eqn} (cf. \cite{DG95}).
	\begin{conj}
		\label{DG conj}
		For fixed coprime integers $ A,B,C \in \Z \setminus \{0\}$, the generalized Fermat equation~\eqref{generalized Fermat eqn} has only finitely many non-trivial coprime integer solutions.
	\end{conj}
	
	In \cite{DG95}, Darmon and Granville proved Conjecture~\ref{DG conj} for fixed $p,q,r$. More precisely,
	\begin{thm}  \rm(\cite[Theorem 2]{DG95})
		For fixed integers $ A,B,C \in \Z \setminus \{0\}$ and fixed $ p,q,r \geq 2$ with $\frac{1}{p} +\frac{1}{q}+ \frac{1}{r} <1$, the generalized Fermat equation~\eqref{generalized Fermat eqn} has only finitely many non-trivial coprime integer solutions.	
	\end{thm}
	
	\subsection{Signature \texorpdfstring{$(p,p,p)$}{(p,p,p)} case}
	Throughout this article, $K$ denotes a totally real number field, $\mcO_K$ denotes its ring of integers, $P:=\mrm{Spec}(\mcO_K)$, $r,n \in \N$, and $p$ denotes a rational prime. Freitas and Siksek~\cite{FS15} first studied asymptotic solutions of the Fermat equation of signature $(p,p,p)$, namely $x^p+y^p+z^p=0$ over $K$ (cf. Definition~\ref{asymptotic solution}). Deconinck~\cite{D16} extended their work to the equation $Ax^p+By^p+Cz^p=0$, where $A,B,C \in \mcO_K$ and $ABC$ is odd (in the sense that $\mfP \nmid ABC$ for every $\mfP \in P$ with $\mfP\mid2$). In \cite{SS18} (resp. \cite{KO20}), \c{S}eng$\ddot{u}$n and Siksek (resp. Kara and Ozman) extended the work of \cite{FS15} (resp. \cite{D16}) to number fields under two standard conjectures (cf. \cite[Conjectures 2.2 and 2.3]{KO20}).
	
	In~\cite{R97}, Ribet proved that the equation $x^p+2^ry^p+z^p=0$ with $2\leq r <p$ has no non-trivial coprime integer solutions.
	In \cite{KS24 Diop1}, Kumar and Sahoo studied the asymptotic solution of the equation $x^p+2^ry^p+z^p=0$ over $K$, for $r \in \N$. 
	
In this article, we study the asymptotic solutions of the generalized Fermat equation of signature $(p,p,p)$, i.e.,
$Ax^p+By^p+Cz^p=0$ over $K$, where $A,B,C\in\mathcal{O}_K\setminus\{0\}$. Deconinck~\cite{D16} studied the asymptotic solutions of this equation under the assumption that $ABC$ is odd. In contrast, we impose no parity restriction on $ABC$. Our principal results are as follows.

\begin{itemize}
	\item In Theorem~\ref{main result1 for Ax^p+By^p+Cz^p=0}, we prove that, for a certain class of fields $K$, the equation
	$Ax^p+By^p+Cz^p=0$
	has no asymptotic solutions in $W_K$ (cf. Definition~\ref{def for W_K} for the definition of $W_K$).
	\item In Theorem~\ref{main result2 for Ax^p+By^p+Cz^p=0}, we prove that, for a certain class of fields $K$, the equation
	$
	Ax^p+By^p+Cz^p=0$
	has no asymptotic solutions in $K^3$, provided that
	$
	A\pm B\pm C\neq 0, \
	\max\{v_{\mfP}(A),v_{\mfP}(BC)\}\leq 4v_{\mfP}(2),$
	and
	$
	v_{\mfP}(ABC)\equiv 0
	\quad\text{ or }
	2v_{\mfP}(2)\pmod{3}$
	for the prime $\mfP\in U_K$ supplied by the $S_K'$-unit condition (cf. \S\ref{notations section for Ax^p+By^p+Cz^p=0} for the definition of $U_K$). Moreover, when $A,B,C\in\mathbb{Z}\setminus\{0\}$, we establish Theorem~\ref{main result2 for Ax^p+By^p+Cz^p=0} without the assumption that $A\pm B\pm C\neq 0$ (cf. Proposition~\ref{relax the assumption}).
	
	\item In Theorem~\ref{thm for density}, we prove that the relevant $S_K$-unit equation has only irrelevant solutions for a density-one set of square-free integers $d\geq2$, where $K=\mathbb{Q}(\sqrt d)$. For coefficients satisfying the additional local hypotheses in Theorem~\ref{main result1 for Ax^p+By^p+Cz^p=0}, respectively Theorem~\ref{main result2 for Ax^p+By^p+Cz^p=0}, this gives the corresponding asymptotic nonexistence result.
	
\end{itemize}
	
	The proofs of Theorems~\ref{main result1 for Ax^p+By^p+Cz^p=0} and~\ref{main result2 for Ax^p+By^p+Cz^p=0} rely on explicit bounds for the solutions of the $S_K^\prime$-unit equation~\eqref{S_K-unit solution} (cf. \S\ref{notations section for Ax^p+By^p+Cz^p=0} for the definition of $S_K^\prime$).
	Furthermore, we present several local criteria for applying Theorems~\ref{main result1 for Ax^p+By^p+Cz^p=0} and~\ref{main result2 for Ax^p+By^p+Cz^p=0}. We use the modular approach to prove these theorems.
	The following are some crucial steps in the modular approach:
	\begin{steps}
		\item For any non-trivial solution $(a, b, c)\in K^3$ to the Diophantine equation $Ax^p+By^p+Cz^p=0$, we attach a \textbf{Frey elliptic curve $E/K$}.
		\item We prove the \textbf{modularity} of $E$ for $p \gg 0$, show that $E$ has \textbf{semistable reduction} at every prime $\mfq \in P$ with $\mfq\mid p$, and prove that the mod-$p$ Galois representation $\bar{\rho}_{E,p}$ is \textbf{irreducible} for $p \gg 0$.
		
		\item 
		Using \textbf{level lowering results} of $\bar{\rho}_{E,p}$, we have $\bar{\rho}_{E,p} \sim \bar{\rho}_{f,p}$, for some Hilbert modular newform defined over $K$ of parallel weight $2$ with rational eigenvalues of lower level.
		\item We show that none of the finitely many Hilbert modular newforms occurring in Step~3 can correspond to $\bar{\rho}_{E,p}$, obtaining a \textbf{contradiction}.
	\end{steps}
	
	\subsection{Signature \texorpdfstring{$(p,p,2)$}{(p,p,2)} case}
	In~\cite{DM97}, Darmon and Merel demonstrated that the equation $x^n+y^n=z^2$ with $n \geq 3$ does not have any non-trivial coprime integer solutions. 
	In~\cite{IKO20}, I\c{s}ik, Kara, and Ozman studied the asymptotic solution of the equation $x^p+y^p=z^2$ of a certain type over $K$ when the narrow class number $h_K^+=1$ and the residual degree $\f(\mfP, 2)=1$ for some prime $\mathfrak{P} \in P$. 
	In ~\cite{KS24 Diop1}, Kumar and Sahoo relaxed the assumptions made in ~\cite{IKO20} and demonstrated that $x^p+y^p=z^2$ has no asymptotic solution in $W_K^\prime$ (cf. \cite[Definition 5.2]{KS24 Diop1} for $W_K^\prime$). 
	In ~\cite{M22}, Mocanu generalized \cite[Theorem 1.1]{IKO20} by replacing the assumption $h_K^+=1$ in ~\cite{IKO20} with $2 \nmid h_K^+$.
	
	In~\cite{I03}, Ivorra examined the integer solutions of equations $x^2=y^p+2^rz^p$ and $2x^2=y^p+2^rz^p$ for $0 \leq r < p$ with primes $p \geq 5$. In ~\cite{S03}, Siksek demonstrated that the only non-trivial coprime integer solutions to the equation $x^2=y^p+2^rz^p$ with exponent $p \geq 5$ and $r \geq 2$ are when $r=3$, $\ x=\pm3,\ y=z=1$. In \cite{BS04}, Bennett and Skinner performed a study on the integer solutions of the generalized Fermat equation~\eqref{generalized Fermat eqn} of signature $(n,n,2)$, i.e., $Ax^n+By^n=Cz^2$, where $A,B,C \in \Z \setminus \{0\}$. In \cite{KS Diop2}, Kumar and Sahoo extended the work of \cite{M22} by investigating the asymptotic solutions of the equation $x^2=By^p+Cz^p$ over $K$, where $B$ is an odd integer and $C$ is either an odd integer or $2^r$ for some $r \in \N$. In \cite{KS Diop2}, the authors also studied the asymptotic solution of the equation $2x^2=By^p+2^rz^p$ over $K$, where $B$ is an odd integer and $r \in \N$.
	
	\subsection{Signature \texorpdfstring{$(p,p,3)$}{(p,p,3)} case}
	In~\cite{DM97}, Darmon and Merel demonstrated that the equation $x^n+y^n=z^3$ with $n \geq 3$ does not have any non-trivial coprime integer solutions. In \cite{BVY04}, Bennett, Vatsal and Yadzdani conducted a study on integer solutions of the generalized Fermat equation~\eqref{generalized Fermat eqn} of signature $(n,n,3)$, i.e., $Ax^n+By^n=Cz^3$, where $A,B,C \in \Z \setminus \{0\}$.
	In \cite{M22}, Mocanu examined the asymptotic solution of the equation $x^p+y^p=z^3$ of a certain type over $K$. 
	In \cite{IKO23},  I\c{s}ik, Kara, and Ozman studied the asymptotic solution of the equation $x^p+y^p=z^3$ of a certain type over any number field $K$ by assuming two standard conjectures under the condition that $h_K^+=1$. Recently in \cite{KS Diop3}, Sahoo and Kumar examined the asymptotic solution of the generalized Fermat equation~\eqref{generalized Fermat eqn} of signature $(p,p,3)$, i.e., $Ax^p+By^p=Cz^3$ over $K$, where $A,B,C \in \mcO_K \setminus \{0\}$. 

\subsection{Limitations of Deconinck's approach and its improvements}
In this subsection, we discuss the method employed in~\cite{D16} to study the asymptotic solutions of the generalized Fermat equation $
Ax^p+By^p+Cz^p=0$ over $K$ under the assumption that $ABC$ is odd. We then explain how this approach can be refined to remove the parity restriction on $ABC$ and thereby treat arbitrary coefficients $A,B,C\in\mcO_K\setminus\{0\}$.
\begin{itemize}
	\item In \cite{D16}, Deconinck used the modularity method inspired by Freitas and Siksek to study the asymptotic Fermat equation $x^p+y^p+z^p=0$ over $K$. The main result in~\cite{FS15} depends on explicit bounds for the solutions of the $S_K$-unit equation, whereas the main result in~\cite{D16} depends on explicit bounds for the solutions of the $S_K'$-unit equation~\eqref{S_K-unit solution} because of the coefficients $A,B,C$ of the equation $Ax^p+By^p+Cz^p=0$ (cf. \S\ref{notations section for Ax^p+By^p+Cz^p=0} for the definition of $S_K'$).	
	
	\item On the other hand, the proofs of our main results, namely Theorems~\ref{main result1 for Ax^p+By^p+Cz^p=0} and~\ref{main result2 for Ax^p+By^p+Cz^p=0}, also rely on explicit bounds for the solutions of the $S_K'$-unit equation~\eqref{S_K-unit solution}. The principal difficulty is that, when $ABC$ is not odd, the reduction behavior at primes $\mfP\in S_K$ of the Frey elliptic curve considered in~\cite{D16} may differ from that in the odd case.
	
	To overcome this difficulty, we first study the reduction type at each prime $\mfP\in S_K$ of the Frey elliptic curve $E$ defined in~\eqref{Frey curve for x^2=By^p+Cz^p of Type I}, associated with a solution $(a,b,c)\in W_K$, where $W_K$ is a distinguished subset of $\mcO_K^3$ (cf. Definition~\ref{def for W_K}). We prove that
	$
	v_{\mfP}(j_E)<0\ 
	\text{and} \
	p\nmid v_{\mfP}(j_E)
	$ for all sufficiently large primes $p$, at primes satisfying the nonvanishing condition in Lemma~\ref{reduction on T and S}(1).
	More generally, for $(a,b,c)\in\mcO_K^3$ and $\mfP\in U_K$, we prove that
	$
	p\mid \#\bar{\rho}_{E,p}(I_{\mfP})$ 
	or $
	3\mid \#\bar{\rho}_{E,p}(I_{\mfP}),
	$ 
	(cf. Lemma~\ref{reduction on T and S}(2)).
	
	\item Since $ABC$ is odd, Deconinck proved the modularity of the Frey curve $E$ for large primes $p$ (cf. \cite[Corollary 2.3]{D16}). On the other hand, we first prove the modularity of the Frey curve $E$ for large primes $p$ when $(a,b,c) \in W_K$. Then, under the condition $A\pm B \pm C \neq 0$, we prove that $E$ is modular for large primes $p$ when $(a,b,c) \in K^3\setminus \{(0,0,0)\}$ (cf. Theorem~\ref{modularity of Frey curve x^2=By^p+Cz^p over W_K}).
	
	\item Since Theorems~\ref{main result1 for Ax^p+By^p+Cz^p=0} and \ref{main result2 for Ax^p+By^p+Cz^p=0} depends on the solutions of the $S_K'$-unit equation~\eqref{S_K-unit solution}. Using the ideas of \cite{FS15} and \cite{KS24 Diop1}, we can provide the local criteria of $K$ for Theorems~\ref{main result1 for Ax^p+By^p+Cz^p=0} and \ref{main result2 for Ax^p+By^p+Cz^p=0} (cf. Propositions~\ref{loc crit for real quadratic field over W_K},~\ref{loc crit for real quadratic field over K},~\ref{loc crit for odd degree field over W_K} and~\ref{local criteria1 for odd defree over K}).
\end{itemize}    

\subsection{Notations}
\label{section for preliminary}
Throughout this article, we fix the following notations.
\begin{itemize}
	\item Let $K$ denote a totally real number field and $ K^\ast:=K\setminus \{0\}$. 
	\item Let $\mcO_K$, $P$, and $\mfn$ be the ring of integers, $\mrm{Spec}(\mcO_K)$, and an ideal of $\mcO_K$, respectively. Let $\mathbb{P}:=\mrm{Spec}(\Z)$.
	\item  For an elliptic curve $E/K$, let $\Delta_E$ and $j_E$ denote the  discriminant and $j$-invariant of $E$, respectively.
	\item Let $G_K:=\Gal(\overline{K}/K)$ denote the absolute Galois group of $K$.
	\item For any prime ideal $\mfP \in P$, let $I_\mfP$ denote the inertia subgroup of $G_K$ at $\mfP$.
	\item For an elliptic curve $E/K$ and a prime $p \in \mathbb{P}$, let $\bar{\rho}_{E,p} : G_K \rightarrow \mathrm{Aut}(E[p]) \simeq \GL_2(\F_p)$ be  the mod-$p$ Galois representation of $G_K$, induced by the action of $G_K$ on the $p$-torsion points $E[p]$ of $E$. 
	\item For any Hilbert modular newform $f$ over $K$ of weight $k$, level $\mfn$ with coefficient field $\Q_f$ and for any $\omega \in \mrm{Spec}(\mcO_{\Q_f})$, let $\bar{\rho}_{f, \omega}: G_K \rightarrow \GL_2(\F_\omega)$ be the residual Galois representation attached to $f, \omega$.
\end{itemize}

\subsection{Structure of the article:}
This article is organized as follows. In \S\ref{notations section for Ax^p+By^p+Cz^p=0}, we state the main results, i.e., Theorems~\ref{main result1 for Ax^p+By^p+Cz^p=0}, ~\ref{main result2 for Ax^p+By^p+Cz^p=0} and~\ref{thm for density} for the asymptotic solutions of the equation $Ax^p+By^p+Cz^p=0$ over $K$. In \S\ref{steps to prove main results}, we prove Theorems~\ref{main result1 for Ax^p+By^p+Cz^p=0} and~\ref{main result2 for Ax^p+By^p+Cz^p=0}. In \S\ref{section for loc criteria}, we provide several local criteria of $K$ such that the equation $Ax^p+By^p+Cz^p=0$ has no asymptotic solution in $W_K$ and $K^3$. Finally, in \S\ref{Proof of the thm for density}, we prove Theorem~\ref{thm for density}.
%

\section{Solutions of the Diophantine equation \texorpdfstring{$Ax^p+By^p+Cz^p=0$}{Axp+Byp+Czp=0} over \texorpdfstring{$K$}{K}}
\label{notations section for Ax^p+By^p+Cz^p=0} 	
In this section, we study the solutions of the following Diophantine equation
\begin{equation}
	\label{Ax^p+By^p+Cz^p=0}
	Ax^p+By^p+Cz^p=0
\end{equation} 
with prime exponent $p\geq 3$ and $A,B,C \in \mcO_K\setminus \{0\}$. 
Let $S_K:= \{ \mfP \in P :\ \mfP|2\}$ and $S_K^{\prime}:= \{ \mfP \in P :\ \mfP|2ABC \}$.

\begin{dfn}[Trivial solution]
	A solution $(a, b, c)\in K^3$ to the equation \eqref{Ax^p+By^p+Cz^p=0} is said to be trivial if $abc=0$,
		and non-trivial otherwise.
	We say $(a, b, c)\in \mcO_K^3$ is primitive if $a\mcO_K+b\mcO_K+c\mcO_K=\mcO_K$.
\end{dfn}  
\begin{dfn}
	\label{def for W_K}
	Let $W_K$ be the set of all non-trivial primitive solutions $(a, b, c)\in  \mcO_K^3$ to equation~\eqref{Ax^p+By^p+Cz^p=0} with $\mfP |abc$ for all $\mfP \in S_K$.
\end{dfn}

\begin{dfn}
	\label{asymptotic solution}
	We say a Diophantine equation $Ax^p+By^p+Cz^p=0$ of prime exponent $p$ has no asymptotic solution in a set $S \subseteq K^3$, if there exists a constant $V_{K,A,B,C}>0$ (depending on $K,A,B,C$) such that for primes $p>V_{K,A,B,C}$, the equation $Ax^p+By^p+Cz^p=0$ has no non-trivial solution in $S$.
\end{dfn}

\begin{remark}
	\label{remark for W_K}
	Let $\mfP \in S_K$. If $(a, b, c)\in W_K$ is a solution to equation~\eqref{Ax^p+By^p+Cz^p=0} with exponent $p > \max\{v_\mfP (A), v_\mfP (B), v_\mfP (C) \}$, then $\mfP$ divides exactly one of $a$, $b$ and $c$. Otherwise, let $\mfP$ divide both $a$ and $b$. Then $\mfP^p | Aa^p+Bb^p=-Cc^p$. Since $p > v_\mfP (C)$, $\mfP |c$, which is not possible because $(a, b, c)$ is primitive. Similarly, the other cases i.e.,  $\mfP$ divide both $b$ and $c$, and $\mfP$ divide both $a$ and $c$, are not possible.
\end{remark}

\subsection{Main results}
\label{section for main result of x^2=By^p+2^rz^p} 
For any set $S \subseteq P$, let $\mcO_{S}:=\{\alpha \in K : v_\mfP(\alpha)\geq 0 \text{ for all } \mfP \in P \setminus S\}$ be the ring of $S$-integers in $K$ and $\mcO_{S}^*$ be the $S$-units of $\mcO_{S}$.   
For $\mfP\in S_K$ and $A,B,C\in\mcO_K\setminus\{0\}$, put
$	\alpha_{\mfP}:=8v_{\mfP}(2)+v_{\mfP}(BCA^{-2}),\ 
	\beta_{\mfP}:=8v_{\mfP}(2)+v_{\mfP}(ACB^{-2}),\
	\gamma_{\mfP}:=8v_{\mfP}(2)+v_{\mfP}(ABC^{-2}).
$
We now show that equation~\eqref{Ax^p+By^p+Cz^p=0} has no asymptotic solution in $W_K$. More precisely,
\begin{thm}
	\label{main result1 for Ax^p+By^p+Cz^p=0}
	Let $K$ be a totally real number field. Let $A,B,C \in \mcO_K\setminus \{0\}$ and $S_K^{\prime}:= \{ \mfP \in P :\ \mfP|2ABC \}$. Suppose that, for every solution $(\lambda, \mu)$ to the $S_K^\prime$-unit equation
	\begin{equation}
		\label{S_K-unit solution}
		\lambda+\mu=1, \ \lambda, \mu \in \mcO_{S_K^\prime}^\ast,
	\end{equation}
	there exists some $\mfP\in S_K$ that satisfies
	\begin{equation}
		\label{assumption for main result x^p+y^p=2^rz^p}
		\max \left\{|v_\mfP(\lambda)|,|v_\mfP(\mu)| \right\}\leq 4v_\mfP(2).
	\end{equation}
	If $\alpha_\mfP\beta_\mfP\gamma_\mfP\ne0$, then the generalized Fermat equation $Ax^p+By^p+Cz^p=0$ has no asymptotic solution in $W_K$.
\end{thm}

\begin{remark}
	In \cite{S14}, Siegel proved that for any finite set $S \subseteq P$, the $S$-unit equation has only a finite number of solutions over any number field $K$, and hence the $S_K^\prime$-unit equation~\eqref{S_K-unit solution} has only a finite number of solutions over $K$. Moreover, these solutions are effectively computable (cf. \cite{AKMRVW21}).
\end{remark}

We write $(ES)$ for ``either $[K: \Q]$ is odd or Conjecture~\ref{ES conj} holds for $K$." Let $U_K:=\{ \mfP \in S_K: 3 \nmid v_\mfP(2) \}$. We now show that equation~\eqref{Ax^p+By^p+Cz^p=0} has no asymptotic solution in $K^3$. More precisely,
\begin{thm}
	\label{main result2 for Ax^p+By^p+Cz^p=0}
	Let $K$ be a totally real number field satisfying the condition $(ES)$. Let $A,B,C \in \mcO_K\setminus \{0\}$ satisfy $A\pm B\pm C\ne0$, and put $S_K^{\prime}:= \{ \mfP \in P :\ \mfP|2ABC \}$.
	Suppose that for every solution $(\lambda, \mu)$ to the $S_K^\prime$-unit equation~\eqref{S_K-unit solution},
	there exists some $\mfP\in U_K$ that satisfies
	\begin{equation}
		\label{assumption for main result2 for Ax^p+By^p+Cz^p=0}
		\max \left\{|v_\mfP(\lambda)|,|v_\mfP(\mu)| \right\}\leq 4v_\mfP(2),
		\qquad v_\mfP(\lambda\mu)\equiv v_\mfP(2) \pmod 3,
	\end{equation}
If 	$
\max\{v_\mfP(A),v_\mfP(BC)\}\leq4v_\mfP(2), \
v_\mfP(ABC)\equiv0\ \text{or}\ 2v_\mfP(2)\pmod3$, then the equation $Ax^p+By^p+Cz^p=0$ has no asymptotic solution in $K^3$.
\end{thm}
The local congruence in Theorem~\ref{main result2 for Ax^p+By^p+Cz^p=0} already implies $\alpha_\mfP\beta_\mfP\gamma_\mfP\ne0$. Indeed, writing $e=v_\mfP(2)$ and $T=v_\mfP(ABC)$, we have
\[
	\alpha_\mfP\equiv\beta_\mfP\equiv\gamma_\mfP
	\equiv 2e+T\equiv e\ \text{or}\ 2e\pmod3,
\]
since $\alpha_\mfP=8e+T-3v_\mfP(A)$, and similarly for $\beta_\mfP$ and $\gamma_\mfP$. As $\mfP\in U_K$, we have $3\nmid e$, so no additional nonvanishing hypothesis is needed.
\begin{remark}
	If $A,B,C \in \Z \setminus \{0\}$, then we prove Theorem~\ref{main result2 for Ax^p+By^p+Cz^p=0} without using the condition $A\pm B \pm C \neq 0$ (cf. Proposition~\ref{relax the assumption}).
\end{remark}
We call $(2,-1)$, $(-1,2)$, and $(\frac{1}{2},\frac{1}{2})$ the irrelevant solutions of the $S_K^\prime$-unit equation~\eqref{S_K-unit solution}; all other solutions are called relevant.
The following is an immediate corollary of Theorems~\ref{main result1 for Ax^p+By^p+Cz^p=0} and~\ref{main result2 for Ax^p+By^p+Cz^p=0}.
\begin{cor}
	\label{irrelevant cor}
	Let $K,A,B,C$ be as above, and suppose that the $S_K^\prime$-unit equation~\eqref{S_K-unit solution} has only irrelevant solutions. If each irrelevant solution admits a prime satisfying all the corresponding hypotheses of Theorem~\ref{main result1 for Ax^p+By^p+Cz^p=0} (resp. Theorem~\ref{main result2 for Ax^p+By^p+Cz^p=0}), then the conclusion of that theorem holds over $K$.
\end{cor}
\subsection{Density computations}
Throughout this subsection, $K=\Q(\sqrt d)$ is a real quadratic field. We calculate the density of the square-free integers $d\geq2$ for which the $S_K$-unit equation has only irrelevant solutions. First, we define relative density.
Let $$\N^{\text{sf}}:=\{d\in \Z_{\geq 2} : d \text{ is a square-free integer}\}.$$
The relative density of $S \subseteq \N^{\text{sf}}$ is defined as follows.
\begin{dfn}
	For $S \subseteq \N^{\text{sf}}$, the relative density of $S$ is defined by 
	$$ \delta_{\text{rel}}(S):= \lim_{x \to \infty}\frac{\# \{d \in S : d \leq x\}}{\# \{d \in \N^{\text{sf}} : d \leq x\}},$$
	if the limit exists.
\end{dfn}
We now state the main result of this subsection.
\begin{thm}
	\label{thm for density}
	For $d\in\N^{\mathrm{sf}}$, put $K_d:=\Q(\sqrt d)$ and let $S_{K_d}$ be the set of primes of $K_d$ above $2$. Let
	\[
		U:=\bigl\{d\in\N^{\mathrm{sf}}:\lambda+\mu=1,
		\ \lambda,\mu\in\mcO_{S_{K_d}}^*,
		\text{ has only irrelevant solutions}\bigr\}.
	\]
	Then $\delta_{\mathrm{rel}}(U)=1$. In particular, let $d\in U$, put $K=K_d$, and suppose that $A,B,C\in\mcO_K\setminus\{0\}$ satisfy $S_K'=S_K$. Then Theorem~\ref{main result1 for Ax^p+By^p+Cz^p=0} applies if there exists $\mfP\in S_K$ such that $\alpha_\mfP\beta_\mfP\gamma_\mfP\ne0$. Likewise, Theorem~\ref{main result2 for Ax^p+By^p+Cz^p=0} applies if $K$ satisfies $(ES)$, $A\pm B\pm C\ne0$, and there is a prime $\mfP\in U_K$ such that
	\[
		\max\{v_\mfP(A),v_\mfP(BC)\}\leq4v_\mfP(2),
		\qquad
		v_\mfP(ABC)\equiv0\ \text{or}\ 2v_\mfP(2)\pmod3.
	\]
\end{thm}
Thus, for almost all real quadratic fields, the $S_K$-unit obstruction required in the main theorems is satisfied. The conclusions of those theorems then follow whenever their respective coefficient hypotheses also hold.
We provide the proof of Theorem~\ref{thm for density} in \S\ref{Proof of the thm for density}.

\section{Steps to prove Theorems~\ref{main result1 for Ax^p+By^p+Cz^p=0} and~\ref{main result2 for Ax^p+By^p+Cz^p=0}}
\label{steps to prove main results}
For any non-trivial solution $(a, b, c)\in K^3$ to equation \eqref{Ax^p+By^p+Cz^p=0} with exponent $p$, consider the Frey elliptic curve 
\begin{equation}
	\label{Frey curve for x^2=By^p+Cz^p of Type I}
	E:=E_{a,b,c} : y^2= x(x-Aa^p)(x+Bb^p),
\end{equation}
where $c_4=2^4(A^2a^{2p}-BCb^pc^p)= 2^4(B^2b^{2p}-ACa^pc^p)=2^4(C^2c^{2p}-ABa^pb^p),\\ \Delta_E=2^4A^2B^2C^2(abc)^{2p}$ and $ j_E=2^{8} \frac{(A^2a^{2p}-BCb^pc^p)^3}{A^2B^2C^2(abc)^{2p}}$. 
\subsection{Modularity of the Frey curve}
In this subsection, we prove the modularity of the Frey curve $E:=E_{a,b,c}$ defined in~\eqref{Frey curve for x^2=By^p+Cz^p of Type I} when $(a,b,c)\in W_K$ (resp. $K^3$) for large primes $p$. First, we recall a modularity result from~\cite{FLHS15}.
\begin{thm}\rm(\cite[Theorem 5]{FLHS15})
	\label{modularity result of elliptic curve over totally real}
	Let $K$ be a totally real number field. Then, up to isomorphism over $\bar{K}$, there exist only finitely many elliptic curves over $K$ that are not modular. 
\end{thm}
\begin{thm}
	\label{modularity of Frey curve x^2=By^p+Cz^p over W_K}
	Let $K$ be a totally real number field.  Then, there exists a constant $D=D_{K,A,B,C}$ (depending on $K,A,B,C$) such that for any non-trivial solution $(a,b,c)\in W_K$ (resp. $K^3$) to equation \eqref{Ax^p+By^p+Cz^p=0} (resp. equation \eqref{Ax^p+By^p+Cz^p=0} with $A\pm B\pm C \neq 0$) with exponent $p >D$, the Frey elliptic curve $E:=E_{a,b,c}$ in~\eqref{Frey curve for x^2=By^p+Cz^p of Type I} is modular.  
\end{thm}

\begin{proof}
	We prove this theorem in two cases.
	\begin{itemize}
		\item Let $(a,b,c)\in W_K$. This case requires substantial new techniques, in contrast to \cite[Corollary 2.3]{D16}, to address the complications arising from the evenness of $ABC$. 
		By Theorem~\ref{modularity result of elliptic curve over totally real}, there exist only finitely many elliptic curves over $K$ up to $\bar{K}$ isomorphism which are not modular. Consider $j_1,\ldots,j_s \in K$ as the $j$-invariants of those elliptic curves.
		Subsequently, the $j$-invariant of the Frey elliptic curve $E$ is given by $ j_E=2^{8} \frac{(A^2a^{2p}-BCb^pc^p)^3}{A^2B^2C^2(abc)^{2p}}=2^{8}\frac{(\mu^2-\mu+1)^3}{\mu^2(\mu-1)^2}$, where $\mu= \frac{-Bb^p}{Aa^p}$.
		For each $i=1,2,\ldots,s$, the equation $j_E=j_i$ has at most six solutions in $K$. So, there exists 
		$\mu_1, \mu_2, ..., \mu_t \in K$ with $t\leq 6s$ such that $E$ is modular for all $\mu \notin\{\mu_1, \mu_2, ..., \mu_t\}$.
		If $\mu=\mu_k$ for some $k\in\{1,2,\ldots,t\}$, then
		\[
			\left(\frac{b}{a}\right)^p=-\frac{A\mu_k}{B},
			\qquad
			\left(\frac{c}{a}\right)^p=\frac{A(\mu_k-1)}{C}.
		\]
		For fixed right-hand sides, only finitely many prime exponents $p$ can occur unless both right-hand sides are roots of unity. In the exceptional case, $b/a$ and $c/a$ are roots of unity. Since $K$ is totally real, $b/a=\pm1$ and $c/a=\pm1$, so $b=\pm a$ and $c=\pm a$. For any fixed $\mfP\in S_K$ and all sufficiently large $p$, this contradicts Remark~\ref{remark for W_K}, because $\mfP$ must divide exactly one of $a,b,c$.
		
		\item Let $(a,b,c) \in K^3$ and assume $A\pm B\pm C \neq 0$. This case follows from \cite[Corollary 2.3]{D16} with minor modifications. As above, only finitely many prime exponents can occur unless both relevant right-hand sides are roots of unity. In that exceptional case, $b/a=\pm1$ and $c/a=\pm1$, which gives a solution $(1,\pm1,\pm1)$ of~\eqref{Ax^p+By^p+Cz^p=0}. It follows that $A\pm B\pm C=0$, a contradiction.
	\end{itemize}
	There are therefore only finitely many exceptional prime exponents. Taking $D$ larger than these exponents and the coefficient bounds above completes the proof.
\end{proof}

\subsection{Irreducibility of the \texorpdfstring{mod-$p$}{mod-p} Galois representations attached to elliptic curves}
Let $E/K$ be an elliptic curve. For any rational prime $p$, let $\mfq \in P$ be a prime ideal of $\mcO_K$ lying above $p$. In~\cite{FS15 Irred}, Freitas and Siksek gave a criterion for the irreducibility of the mod-$p$ Galois representation $\bar{\rho}_{E,p}$ for large primes $p$. More precisely,

\begin{thm} \rm(\cite[Theorem 2]{FS15 Irred})
	\label{irreducibility of mod $P$ representation}
	Let $K$ be a totally real Galois field. Then there exists an effective constant $C_K$ such that, if $p>C_K$ is a prime and $E/K$ is an elliptic curve that is semistable at every $\mfq\mid p$, then $\bar{\rho}_{E,p}$ is irreducible.
\end{thm}	

\subsection{Level lowering results}
For any $\mfq \in P$, let $\Delta_\mfq$ be the minimal discriminant of $E$ at $\mfq$. Let 

\begin{equation}
	\label{conductor of elliptic curve}
	\mfm_p:= \prod_{ p|v_\mfq(\Delta_\mfq), \ \mfq ||\mfn} \mfq \text{ and } \mfn_p:=\frac{\mfn}{\mfm_p}.
\end{equation}

In~\cite{FS15}, Freitas and Siksek gave a level-lowering result using the work of Fujiwara~\cite{F06}, Jarvis~\cite{J04}, and Rajaei~\cite{R01}. More precisely,
\begin{thm} \rm(\cite[Theorem 7]{FS15})
	\label{level lowering of mod $p$ repr}
	Let $E$ be an elliptic curve over $K$ of conductor
	$\mfn$. Let $p$ be a rational prime. Suppose that the following conditions hold:
	\begin{enumerate}
		\item  For $p \geq 5$, the ramification index $e(\mfq /p) < p-1$ for all $\mfq |p$, and $\Q(\zeta_p)^+ \nsubseteq K$;
		\item $E/K$ is modular;
		\item $\bar{\rho}_{E,p}$ is irreducible;
		\item $E$ is semistable at every $\mfq\mid p$, and $p\mid v_\mfq(\Delta_\mfq)$ for every $\mfq\mid p$.
	\end{enumerate}
	Then there exists a Hilbert modular newform $f$ over $K$ of parallel weight $2$, level $\mfn_p$, and some prime $\lambda$ of $\mcO_{Q_f}$ such that $\lambda | p$ and $\bar{\rho}_{E,p} \sim \bar{\rho}_{f,\lambda}$.
\end{thm}	

\subsection{Eichler-Shimura}
Let $E/K$ be an elliptic curve of conductor $\mfn$.
We now state a conjecture, which is an extension of the Eichler-Shimura theorem over $\Q$. 
\begin{conj}[Eichler-Shimura]
	\label{ES conj}
	Let $f$ be a Hilbert modular newform over $K$ of parallel weight $2$, level $\mfn$, and with coefficient field $\Q_f= \Q$. Then, there exists an elliptic curve $E_f /K$ with conductor $\mfn$ having the same $L$-function as $f$.
\end{conj}
In~\cite[Theorem 7.7]{D04}, Darmon showed that Conjecture~\ref{ES conj} holds over $K$, if  either 
$[K: \Q] $ is odd or there exists some prime ideal $\mfq \in P$ such that $v_\mfq(\mfn) = 1$. 
In \cite{FS15}, Freitas and Siksek provided a partial answer to Conjecture~\ref{ES conj} in terms of mod-$p$ Galois representations attached to $E$. More precisely,

\begin{thm} \rm(\cite[Corollary 2.2]{FS15})
	\label{FS partial result of E-S conj}
	Let $E$ be an elliptic curve over $K$ and $p$ be an odd prime.
	Suppose that $\bar{\rho}_{E,p}$ is irreducible and $\bar{\rho}_{E,p} \sim \bar{\rho}_{f,p}$ for some Hilbert modular newform $f$ over $K$ of parallel weight $2$ and level
	$\mfn$ with rational eigenvalues.
	Let $\mfq \nmid p$ be a prime of $K$ such that
	\begin{enumerate}
		\item $E$ has potentially multiplicative reduction at $\mfq$ (i.e., $v_\mfq(j_E)<0$);
		\item $p| \# \bar{\rho}_{E,p}(I_\mfq)$;
		\item  $p \nmid \left(\mathrm{Norm}(K/\Q)(\mfq) \pm 1\right)$.
	\end{enumerate}
	Then there exists an elliptic curve $E_f /K$ of conductor $\mfn$ having the same $L$-function as $f$.
\end{thm}
\subsection{Conductor of the Frey curve}
Any non-trivial solution $(a,b,c)\in K^3$ of the generalized Fermat equation $Ax^p+By^p+Cz^p=0$ can be rescaled so that $(a,b,c)\in\mcO_K^3$. Since the class number $h_K$ need not be $1$, we cannot in general assume that $a,b,c$ are coprime. We can, however, arrange that the ideal $a\mcO_K+b\mcO_K+c\mcO_K$ belongs to a fixed finite set $H$ of prime ideals of $\mcO_K$. Thus $a,b,c$ are coprime away from $H\cup\{\mfp\in P:\mfp\mid ABC\}$. We now explain this in more detail.

For any nonzero fractional ideal $\mfa$ of $\mcO_K$, the image of $\mfa$ in the class group of $K$ is denoted by $[\mfa]$.
For any non-trivial solution $(a,b,c) \in K^3$ to equation $Ax^p+By^p+Cz^p=0$, let $G_{a,b,c}:=a\mcO_K+  
b\mcO_K+c\mcO_K$. Let $[a,b,c] :=[G_{a,b,c}]$. 
Let the class number be $h_K=h$, and let $\mfc_1,\mfc_2,\ldots,\mfc_h$ be the distinct ideal classes of $K$. Every ideal class contains infinitely many prime ideals of $\mcO_K$ (see \cite[\S1.8]{FS15}). Let $H$ be a set of prime ideals $\mfm_i\nmid2ABC$, one representing each ideal class $\mfc_i$. Recall that $S_K^{\prime}:=\{\mfP\in P:\mfP\mid2ABC\}$. Clearly, $S_K^{\prime}\cap H=\varnothing$.
%

The following lemma, similar to \cite[Lemma 3.2]{FS15} and \cite[Lemma 3.1]{D16}, states that we can scale any solution $(a,b,c)\in K^3$ of~\eqref{Ax^p+By^p+Cz^p=0} to an integral solution with $G_{a,b,c}\in H$. More precisely,
\begin{lem}
	\label{integrality of solution}
	Let $(a,b,c) \in K^3$ be a non-trivial solution of equation~\eqref{Ax^p+By^p+Cz^p=0}. Then, there exists a non-trivial solution $(a',b',c') \in \mcO_K^3$ of equation~\eqref{Ax^p+By^p+Cz^p=0} such that the following conditions hold.
	\begin{enumerate}
		\item $(a^\prime,b^\prime,c^\prime)=(\xi a,\xi b,\xi c)$ for some $\xi\in K^\ast$.
		\item $G_{a^\prime,b^\prime,c^\prime}=\mfm$ for some $\mfm \in H$.
		\item $[a^\prime,b^\prime,c^\prime]=[a,b,c]$.
	\end{enumerate}
\end{lem}


\begin{remark}
	\label{remark for K}
		Let $\mfP \in S_K$, and let $(a,b,c) \in K^3$ be a non-trivial solution of~\eqref{Ax^p+By^p+Cz^p=0} with exponent $p > \max\{v_\mfP(A),v_\mfP(B),v_\mfP(C)\}$. By Lemma~\ref{integrality of solution}, we can rescale $(a,b,c)$ so that $(a,b,c)\in\mcO_K^3$ and $G_{a,b,c}=\mfm$ for some $\mfm\in H$. If $\mfP\mid abc$, then $\mfP$ divides exactly one of $a,b,c$. Indeed, if $\mfP$ divided both $a$ and $b$, then it would also divide $c$ because $p>v_\mfP(C)$. This would imply $\mfP\mid\mfm$, which is impossible because $\mfm\nmid2$. The other cases are identical.
\end{remark}

The following lemma characterizes the reduction of the Frey curve away from the primes in its lowered level and determines its conductor.
\begin{lem}
	\label{reduction away from S}
	Let $(a,b,c) \in W_K$ (resp. $\mcO_K^3$) be a non-trivial solution to equation~\eqref{Ax^p+By^p+Cz^p=0} with prime exponent $p$ (resp. with $G_{a,b,c}=\mfm$ for some $\mfm \in H$), and let $E$ be the Frey curve given in~\eqref{Frey curve for x^2=By^p+Cz^p of Type I}. Put $\Sigma=S_K'$ in the first case and $\Sigma=S_K'\cup\{\mfm\}$ in the second. Then, at every prime $\mfq\in P\setminus\Sigma$, the curve $E$ is minimal and semistable at $\mfq$, and $p\mid v_\mfq(\Delta_E)$. Let $\mfn$ be the conductor of $E$, and let $\mfn_p$ be as in~\eqref{conductor of elliptic curve}. Then
	\begin{equation}
		\label{conductor of E and E' x^2=By^p+Cz^p Type I}
		\mfn=\prod_{\mfP\in\Sigma}\mfP^{r_\mfP}
		\prod_{\substack{\mfq\mid abc\\ \mfq\notin\Sigma}}\mfq,
		\qquad
		\mfn_p=\prod_{\mfP\in\Sigma}\mfP^{r_\mfP^{\prime}},
	\end{equation}
	where $0 \leq r_\mfP^{\prime} \leq r_\mfP $ with $r_\mfP 
	\leq 2+6v_\mfP(2)$ for $\mfP |2$ and $ r_\mfP\leq 2+3v_\mfP(3)$ for $\mfP \nmid 2$.
\end{lem}

\begin{proof}
	The arguments in the proof of this lemma follow directly from \cite[Lemma 3.2]{D16} with minor modifications, and we include them here for completeness.
	Let $\mfq\in P\setminus\Sigma$.
	If $\mfq \nmid \Delta_E$, then $E$ has good reduction at $\mfq$
	and $p | v_\mfq(\Delta_E)=0$. 
	
	Otherwise, suppose $\mfq\mid\Delta_E=2^4A^2B^2C^2(abc)^{2p}$. Since $\mfq\notin\Sigma$, we have $\mfq\nmid2ABC$. If $(a,b,c) \in W_K$, primitivity implies that $\mfq$ divides exactly one of $a,b,c$. In the second case, the same conclusion follows from $G_{a,b,c}=\mfm$ and $\mfq\ne\mfm$. Therefore $\mfq\nmid c_4$, so $E$ is minimal and has multiplicative reduction at $\mfq$.
	Since $v_\mfq(\Delta_E)=2p v_\mfq(abc) $, $p | v_\mfq(\Delta_E)$. 
	By the definition of $\mfn_p$ in~\eqref{conductor of elliptic curve}, we conclude that $\mfq\nmid\mfn_p$ for every $\mfq\notin\Sigma$. Finally, for $\mfP\in\Sigma$, the bounds on $r_\mfP$ follow from \cite[Theorem IV.10.4]{S94}.
\end{proof}

\subsection{Type of reduction with image of inertia}
The following lemma is helpful for the types of reduction of the Frey curve at $\mfP \in P$. 
\begin{lem} \rm(\cite[Lemma 3.4]{FS15})
	\label{criteria for potentially multiplicative reduction}
	Let $E/K$ be an elliptic curve and $p>5$ be a rational prime. For prime $\mfq \in P$ with $\mfq \nmid p$, $E$ has potentially multiplicative reduction at $\mfq$  and $p \nmid v_\mfq(j_E)$ if and only if $p | \# \bar{\rho}_{E,p}(I_\mfq)$.
\end{lem}

The following lemma determines the types of reduction of the Frey curve $E$ at $\mfm$.
\begin{lem}
	\label{Type of reduction at m}
	Let $(a,b,c) \in \mcO_K^3$ be a non-trivial solution to equation~\eqref{Ax^p+By^p+Cz^p=0} with exponent $p>5$ and with $G_{a,b,c}=\mfm$ for some $\mfm \in H$. If $\mfm\nmid p$, then $\#\bar{\rho}_{E,p}(I_\mfm)\mid24$ and hence $p \nmid \#\bar{\rho}_{E,p}(I_\mfm)$.
\end{lem}

\begin{proof}
	Recall that
	$
		\Delta_E=2^4A^2B^2C^2(abc)^{2p}
	$
	and
	$
		c_4=2^4(A^2a^{2p}-BCb^pc^p)
		=2^4(B^2b^{2p}-ACa^pc^p)
		=2^4(C^2c^{2p}-ABa^pb^p).
	$
	Since $\mfm\nmid2ABC$, if $a,b,c$ have unequal valuations at $\mfm$, then
	$
		v_\mfm(c_4)=2p\min\{v_\mfm(a),v_\mfm(b),v_\mfm(c)\}$ and $
		v_\mfm(\Delta_E)=2pv_\mfm(abc).$
	Consequently, $v_\mfm(j_E)<0$ and $p\mid v_\mfm(j_E)$. Lemma~\ref{criteria for potentially multiplicative reduction} therefore gives $p\nmid\#\bar{\rho}_{E,p}(I_\mfm)$. If the three valuations have the common value $t$, then $v_\mfm(c_4)\geq2pt$ and $v_\mfm(\Delta_E)=6pt$, so $v_\mfm(j_E)\geq0$; the same lemma again gives $p\nmid\#\bar{\rho}_{E,p}(I_\mfm)$. In the first case, the inertia image has order at most $2$, while in the potentially good case its order divides $24$. Thus $\#\bar{\rho}_{E,p}(I_\mfm)\mid24$ in either case.
\end{proof}

The following lemma is useful for determining the type of reduction of the Frey curve $E$ at $\mfP \in U_K$.
\begin{lem} \rm(\cite[Lemma 3.6]{FS15})
	\label{3 divides discriminant}
	Let $E/K$ be an elliptic curve. Let $p\geq 3$ be a prime and $\mfP \in S_K$. Suppose that $E$ has potentially good reduction at $\mfP$ (i.e., $v_\mfP(j_E)\geq0$). Then $3 \nmid v_\mfP(\Delta_E)$ if and only if $3\mid\#\bar{\rho}_{E,p}(I_\mfP)$.
\end{lem}

Recall that for any $\mfP \in S_K$, we have $\alpha_{\mfP}
:=8v_{\mfP}(2)+v_{\mfP}(BCA^{-2}),\ 
\beta_{\mfP}:=8v_{\mfP}(2)+v_{\mfP}(ACB^{-2}),\
\gamma_{\mfP}:=8v_{\mfP}(2)+v_{\mfP}(ABC^{-2})$. We will now determine the type of reduction of the Frey curve $E_{a,b,c}$ at $\mfP \in S_K$. 
\begin{lem}
	\label{reduction on T and S}
	Let $\mfP\in S_K$. Let $(a,b,c)\in\mcO_K^3$ be a non-trivial solution to equation~\eqref{Ax^p+By^p+Cz^p=0} with prime exponent
	$
		p > \max \Big\{5, v_\mfP(ABC), |\alpha_{\mfP}|, |\beta_{\mfP}|, |\gamma_{\mfP}| \Big\}.
	$
	Suppose that either $(a,b,c)\in W_K$ or $G_{a,b,c}=\mfm$ for some $\mfm\in H$, and let $E$ be the Frey curve given in~\eqref{Frey curve for x^2=By^p+Cz^p of Type I}.
	\begin{enumerate}
		\item If $(a,b,c)\in W_K$ and $\alpha_\mfP\beta_\mfP\gamma_\mfP\ne0$, then $\ v_\mfP(j_E) < 0$ and $p \nmid v_\mfP(j_E)$, equivalently $p | \#\bar{\rho}_{E,p}(I_\mfP)$.
		\item Assume
		$
			\max\{v_\mfP(A),v_\mfP(BC)\}\leq4v_\mfP(2)
		$
		and $v_\mfP(ABC)\equiv0$ or $2v_\mfP(2)\pmod3$. If $(a,b,c)\in\mcO_K^3$ and $\mfP\in U_K$, then either $p\mid\#\bar{\rho}_{E,p}(I_\mfP)$ or $3\mid\#\bar{\rho}_{E,p}(I_\mfP)$.
	\end{enumerate}
\end{lem} 
\begin{proof}
	Recall that $ j_E=2^{8} \frac{(A^2a^{2p}-BCb^pc^p)^3}{A^2B^2C^2(abc)^{2p}}$ and $\Delta_E=2^4A^2B^2C^2(abc)^{2p}$. 
	\begin{enumerate}
		\item Let $(a,b,c)\in W_K$ and assume $\alpha_\mfP\beta_\mfP\gamma_\mfP\ne0$. By Remark~\ref{remark for W_K}, $\mfP$ divides exactly one of $a$, $b$ and $c$. 
		Suppose $\mfP |a$. Then $\mfP \nmid bc$. 
		Since $p >v_\mfP(ABC) \geq v_\mfP(BC)$, we get $v_\mfP(j_E)=8v_\mfP(2)+3v_\mfP(BC)-2v_\mfP(ABC)-2pv_\mfP(a)=\alpha_\mfP-2pv_\mfP(a)$. Since $p>|\alpha_\mfP|$ and $\alpha_\mfP\ne0$, we have $v_\mfP(j_E)<0$ and $p\nmid v_\mfP(j_E)$. The cases $\mfP\mid b$ and $\mfP\mid c$ are identical, using $\beta_\mfP$ and $\gamma_\mfP$, respectively. Finally, Lemma~\ref{criteria for potentially multiplicative reduction} gives $p\mid\#\bar{\rho}_{E,p}(I_\mfP)$.
		
		\item Let $\mfP\in U_K$, and put $e=v_\mfP(2)$ and $T=v_\mfP(ABC)$. Since $\alpha_\mfP=8e+T-3v_\mfP(A)$, and similarly for $\beta_\mfP$ and $\gamma_\mfP$, the coefficient congruence gives
		\[
			\alpha_\mfP\equiv\beta_\mfP\equiv\gamma_\mfP
			\equiv2e+T\equiv e\ \text{or}\ 2e\pmod3.
		\]
		As $3\nmid e$, it follows that $\alpha_\mfP\beta_\mfP\gamma_\mfP\ne0$. If $\mfP\mid abc$, then, using Remarks~\ref{remark for W_K} and~\ref{remark for K} and the local calculation in (1), we obtain $p\mid\#\bar{\rho}_{E,p}(I_\mfP)$. 
		
		Suppose $\mfP \nmid abc$. If $2v_\mfP(A) < v_\mfP(BC)$, then $v_\mfP(j_E)=  8 v_\mfP(2)+6 v_\mfP(A)-2v_\mfP(ABC)= 8 v_\mfP(2)+4 v_\mfP(A)-2v_\mfP(BC) $. Since $v_\mfP(BC) \leq 4v_\mfP(2)$, it follows that $v_\mfP(j_E) \geq 0$.
		If $2v_\mfP(A) \geq v_\mfP(BC)$, then $v_\mfP(j_E)\geq  8 v_\mfP(2)+3 v_\mfP(BC)-2v_\mfP(ABC)=  8 v_\mfP(2)+ v_\mfP(BC)-2v_\mfP(A)$. Since $v_\mfP(A) \leq 4v_\mfP(2)$, it follows that $v_\mfP(j_E) \geq 0$.
		Here $v_\mfP(\Delta_E)= 4v_\mfP(2)+2v_\mfP(ABC)$. Now, using the hypothesis $v_\mfP(ABC) \equiv 0$ or $2v_\mfP(2) \pmod 3$, we get $v_\mfP(\Delta_E) \equiv v_\mfP(2)$ or $2v_\mfP(2) \pmod 3$. Since $\mfP \in U_K$, we get $3 \nmid v_\mfP(\Delta_E)$. Finally, the proof of $(2)$ follows from Lemma~\ref{3 divides discriminant}.
	\end{enumerate}
\end{proof}

\subsection{Proof of Theorem~\ref{main result1 for Ax^p+By^p+Cz^p=0}.}
The proof of this theorem depends on the following result.	
\begin{thm}
	\label{auxilary result x^2=By^p+Cz^p over W_K}
	Let $K$ be a totally real number field, let $A,B,C\in\mcO_K\setminus\{0\}$, and suppose that there exists $\mfP_0\in S_K$ such that $\alpha_{\mfP_0}\beta_{\mfP_0}\gamma_{\mfP_0}\ne0$. Then there is a constant $V=V_{K,A,B,C}>0$ such that the following holds.
	Let $(a,b,c)\in W_K$ be a solution to equation \eqref{Ax^p+By^p+Cz^p=0} with exponent $p>V$, and let $E$ be the Frey curve given in \eqref{Frey curve for x^2=By^p+Cz^p of Type I}. Then, there exists an elliptic curve $E^\prime/K$ such that:
	\begin{enumerate}
		\item $E^\prime$ has good reduction away from $S_K^{\prime}$ and has full $2$-torsion;
		\item $\bar{\rho}_{E,p} \sim\bar{\rho}_{E^\prime,p}$;
		\item $v_\mfP(j_{E^\prime})<0$ for every $\mfP\in S_K$ such that $\alpha_\mfP\beta_\mfP\gamma_\mfP\ne0$.
	\end{enumerate}
\end{thm}

\begin{proof}[Proof of Theorem~\ref{auxilary result x^2=By^p+Cz^p over W_K}] 
	By Theorem~\ref{modularity of Frey curve x^2=By^p+Cz^p over W_K}, the Frey curve $E$ is modular for primes $p\gg0$. By Lemma~\ref{reduction away from S}, $E$ is semistable away from $S_K^\prime$. Using Theorem~\ref{irreducibility of mod $P$ representation} and replacing $K$ with its Galois closure, we conclude that $\bar{\rho}_{E,p}$ is irreducible for all primes $p \gg 0$.
	Now, by Theorem~\ref{level lowering of mod $p$ repr}, there exist a Hilbert modular newform $f$ over $K$ of parallel weight $2$, level $\mfn_p$ and some prime ideal $\lambda$ of $\mcO_{\Q_f}$ with $\lambda | p$ such that $\bar{\rho}_{E,p} \sim \bar{\rho}_{f,\lambda}$ for $p \gg 0$. 
	By allowing $p$ to be sufficiently large, we can assume $\Q_f=\Q$ (cf. \cite[\S 4]{FS15} for more details).

	Choose $\mfP_0\in S_K$ such that $\alpha_{\mfP_0}\beta_{\mfP_0}\gamma_{\mfP_0}\ne0$. By Lemma~\ref{reduction on T and S}, $E$ has potentially multiplicative reduction at $\mfP_0$ and $p\mid\#\bar{\rho}_{E,p}(I_{\mfP_0})$
	for $p \gg 0$. The existence of an elliptic curve $E_f$ of conductor $\mfn_p$
	then follows from Theorem~\ref{FS partial result of E-S conj} for all $p \gg 0$ (after excluding the primes $p \mid \left( \text{Norm}(K/\Q)(\mfP_0) \pm 1 \right)$). Consequently, $\bar{\rho}_{E,p} \sim \bar{\rho}_{E_f,p}$ for all primes $p>V:=V_{K,A,B,C}$, where $V_{K,A,B,C}$ is the maximum of all the above lower bounds and the finitely many bounds in Lemma~\ref{reduction on T and S} for primes $\mfP\in S_K$ satisfying $\alpha_\mfP\beta_\mfP\gamma_\mfP\ne0$.
	
	\begin{enumerate}
		\item Since the conductor of $E_f$ is $\mfn_p$ given in \eqref{conductor of E and E' x^2=By^p+Cz^p Type I}, $E_f$ has good reduction away from $S_K^\prime$. After enlarging $V$ by an effective amount and by possibly replacing $E_f$ with an isogenous curve, say $E^\prime$, we find that $E^\prime/ K$ has full $2$-torsion. 
		This follows from~\cite[Proposition 15.4.2]{C07} and the fact that $E/K$ has all its points of order $2$ (cf.~\cite[\S 4]{FS15}, \cite[\S3.8]{JS25} for more details). Since $E_f$ is isogenous to $E^\prime$, it follows that $E^\prime$ has good reduction away from $S_K^\prime$.
		
		\item Since $E_f$ is isogenous to $E^\prime$, $\bar{\rho}_{E,p} \sim \bar{\rho}_{E^\prime,p}$ for primes $p>V$.
		\item Let $\mfP\in S_K$ satisfy $\alpha_\mfP\beta_\mfP\gamma_\mfP\ne0$. Lemma~\ref{reduction on T and S} gives $p\mid\#\bar{\rho}_{E,p}(I_\mfP)=\#\bar{\rho}_{E^\prime,p}(I_\mfP)$. Therefore, Lemma~\ref{criteria for potentially multiplicative reduction} gives $v_\mfP(j_{E^\prime})<0$.
	\end{enumerate}
	This completes the proof of the theorem.
\end{proof}
We now prove Theorem~\ref{main result1 for Ax^p+By^p+Cz^p=0}, and its proof is similar to that of ~\cite[Theorem 3.3]{KS24 Diop1}.
\begin{proof}[Proof of Theorem~\ref{main result1 for Ax^p+By^p+Cz^p=0}]
	Let $V=V_{K,A,B,C}$ be as in Theorem~\ref{auxilary result x^2=By^p+Cz^p over W_K}, and let $(a,b,c)\in W_K$ be a solution to equation~\eqref{Ax^p+By^p+Cz^p=0} with exponent $p>V$. By that theorem, there exists an elliptic curve $E^\prime/K$ with full $2$-torsion and good reduction away from $S_K^\prime$.
	Hence $E^\prime$ has a model of the form $E^\prime: Y^2 = (X-e_1)(X-e_2)(X-e_3)$, where $e_1, e_2, e_3$ are distinct and their cross ratio $\lambda= \frac{e_3-e_1}{e_2-e_1} \in \mathbb{P}^1(K)-\{0,1,\infty\}$. Then $E^\prime$ is isomorphic (over $\bar{K}$) to an elliptic curve $E_\lambda$ in the Legendre form:
	$$E_\lambda : y^2 = x(x - 1)(x-\lambda) \text{ for } \lambda \in \mathbb{P}^1(K)-\{0,1,\infty\}$$ 
	with
	\begin{equation}
		\label{j'-invariant of Legendre form}
		j_{E^\prime} = j(E_\lambda)= 2^8\frac{(\lambda^2-\lambda+1)^3}{\lambda^2(1-\lambda)^2}.
	\end{equation}
	Then the action of the symmetric group $S_3$ on $\{e_1,e_2,e_3\}$ can be extended to an action on $\mathbb{P}^1(K)-\{0,1,\infty\}$ via the cross ratio $\lambda= \frac{e_3-e_1}{e_2-e_1}$. Under the action of $S_3$, the orbit
	of $\lambda \in \mathbb{P}^1(K)-\{0,1,\infty\}$ is the set $\left \{\lambda, \frac{1}{\lambda}, 1-\lambda, \frac{1}{1-\lambda}, \frac{\lambda}{\lambda-1}, \frac{\lambda-1}{\lambda}\right \}$,
	which are called the $\lambda$-invariants of $E'$ (cf. \cite[Page 10]{KS24 Diop1}, \cite[\S5]{FS15} for more details).
	
	Since $E^\prime$ has good reduction away from $S_K^\prime$, $j_{E^\prime} \in \mcO_{S_K^\prime}$. Let $\lambda\in K$ be any $\lambda$-invariant of $E'$. By \eqref{j'-invariant of Legendre form}, $\lambda\in K$ satisfies a monic polynomial equation of degree $6$ over $\mcO_{S_K^\prime}$, and hence $\lambda \in \mcO_{S_K^\prime}$. Similarly $\frac{1}{\lambda}$,  $\mu:=1-\lambda$, and $\frac{1}{\mu}$ belongs to $\mcO_{S_K^\prime}$. Therefore, $\lambda, \mu \in \mcO^*_{S_K^\prime}$. Hence $(\lambda, \mu)$ is a solution of the $S_K^\prime$-unit equation~\eqref{S_K-unit solution}. 
	We now rewrite \eqref{j'-invariant of Legendre form} in terms of $\lambda, \ \mu$ as
	\begin{equation}
		\label{j' in terms of lambda and mu}
		j_{E^\prime}= 2^8\frac{(1-\lambda \mu)^3}{(\lambda \mu)^2}.
	\end{equation}
	
	By the hypotheses of Theorem~\ref{main result1 for Ax^p+By^p+Cz^p=0}, there exists some $\mfP\in S_K$ such that
	\[
		\alpha_\mfP\beta_\mfP\gamma_\mfP\ne0,
		\qquad
		s:=\max\{|v_\mfP(\lambda)|,|v_\mfP(\mu)|\}\leq4v_\mfP(2).
	\]
	If $s=0$, then $v_\mfP(\lambda)= v_\mfP(\mu)=0$. This gives $v_\mfP(j_{E^\prime})\geq 8v_\mfP(2)>0$, which contradicts Theorem~\ref{auxilary result x^2=By^p+Cz^p over W_K}$(3)$. 
	On the other hand, let $s>0$. Since $\lambda+\mu=1$, either $v_\mfP(\lambda)=v_\mfP(\mu)=-s$, $v_\mfP(\lambda)=0$ and $v_\mfP(\mu)=s$, or $v_\mfP(\lambda)=s$ and $v_\mfP(\mu)=0$. Thus $v_\mfP(\lambda\mu)=-2s$ or $s$. In every case, $v_\mfP(j_{E^\prime})\geq8v_\mfP(2)-2s\geq0$, contradicting Theorem~\ref{auxilary result x^2=By^p+Cz^p over W_K}(3).
\end{proof}

\subsection{Proof of Theorem~\ref{main result2 for Ax^p+By^p+Cz^p=0}.}
The proof of this theorem depends on the following result.	
\begin{thm}
	\label{auxilary result Ax^p+By^p+Cz^p=0 over K}
	Let $K$ be a totally real number field satisfying $(ES)$, and let $A,B,C\in\mcO_K\setminus\{0\}$ satisfy $A\pm B\pm C\ne0$. Then there is a constant $V=V_{K,A,B,C}>0$ such that the following holds.
	Let $(a,b,c)\in K^3$ be a non-trivial solution to equation \eqref{Ax^p+By^p+Cz^p=0} with exponent $p>V$, and rescale $(a, b, c)$ so that $(a, b, c) \in \mcO_K^3$ and $G_{a,b,c}=\mfm$ for some $\mfm \in H$.
	Let $E$ be the Frey curve given in \eqref{Frey curve for x^2=By^p+Cz^p of Type I}. Then, there exists an elliptic curve $E^\prime/K$ such that:
	\begin{enumerate}
		\item $E^\prime$ has good reduction away from $S_K^{\prime} \cup \{\mfm\}$;
		\item $E^\prime/K$ has full $2$-torsion and $\bar{\rho}_{E,p} \sim\bar{\rho}_{E^\prime,p}$;
		\item $E'$ has potentially good reduction at $\mfm$;
		\item If $\mfP\in U_K$, $\max\{v_\mfP(A), v_\mfP(BC)\} \leq 4v_\mfP(2)$, and $v_\mfP(ABC) \equiv 0$ or $2v_\mfP(2) \pmod 3$, then either $v_\mfP(j_{E^\prime})<0$ or $3 \nmid v_\mfP(j_{E^\prime})$.
	\end{enumerate}
\end{thm}
For $A,B,C\in\Z\setminus\{0\}$, the assumption $A\pm B\pm C\ne0$ may be omitted, by the modularity argument in Proposition~\ref{relax the assumption}.

\begin{proof}
	The proofs of (1) and (2) are similar to the proof of Theorem~\ref{auxilary result x^2=By^p+Cz^p over W_K}. Enlarge $V$ so that the prime bound in Lemma~\ref{reduction on T and S} is satisfied at every $\mfP\in S_K$ and so that no $\mathfrak m\in H$ lies above a rational prime $p>V$. There are only finitely many $K$-isomorphism classes of elliptic curves over $K$ with full $2$-torsion and good reduction away from $S_K'\cup\{\mathfrak m\}$. Since $H$ is finite, after enlarging $V$ once more, we may therefore assume that
	$
		v_{\mathfrak m}(j_{E'})<0$ and $
		p\nmid v_{\mathfrak m}(j_{E'})$
	for every $\mathfrak m\in H$ and every prime $p>V$. By Lemma~\ref{Type of reduction at m},
	$
		p\nmid\#\bar{\rho}_{E,p}(I_{\mathfrak m}).
	$
	As $\bar{\rho}_{E,p}\cong\bar{\rho}_{E',p}$, the same is true for $E'$. If $v_{\mathfrak m}(j_{E'})<0$, Lemma~\ref{criteria for potentially multiplicative reduction} would instead imply that $p\mid\#\bar{\rho}_{E',p}(I_{\mathfrak m})$, a contradiction. Hence $v_{\mathfrak m}(j_{E'})\ge0$, proving (3).

	Let $\mfP\in U_K$ satisfy the coefficient conditions in (4). By Lemma~\ref{reduction on T and S}(2), either $p\mid\#\bar{\rho}_{E,p}(I_\mfP)$ or $3\mid\#\bar{\rho}_{E,p}(I_\mfP)$.
	\begin{itemize}
		\item Suppose $p | \# \bar{\rho}_{E,p}(I_\mfP)= \# \bar{\rho}_{E^\prime,p}(I_\mfP)$. Then by Lemma~\ref{criteria for potentially multiplicative reduction}, we get $v_\mfP(j_{E^\prime})<0$, hence we are done.
		\item Suppose $p\nmid\#\bar{\rho}_{E,p}(I_\mfP)$. By Lemma~\ref{reduction on T and S},
		$
			3\mid\#\bar{\rho}_{E,p}(I_\mfP)
			=\#\bar{\rho}_{E^\prime,p}(I_\mfP).
		$
		If $v_\mfP(j_{E^\prime})<0$, then we are done. Otherwise, $v_\mfP(j_{E^\prime})\geq0$, and Lemma~\ref{3 divides discriminant} gives $3\nmid v_\mfP(\Delta_{E^\prime})$.
		Since $v_\mfP(j_{E^\prime}) = 3v_\mfP(c_4) -v_\mfP(\Delta_{E^\prime})$ and  $3\nmid v_\mfP(\Delta_{E^\prime})$, we have $3\nmid v_\mfP(j_{E^\prime})$.
	\end{itemize}
	This completes the proof of the theorem.
\end{proof}

\begin{proof}[Proof of Theorem~\ref{main result2 for Ax^p+By^p+Cz^p=0}]
	Let $V=V_{K,A,B,C}$ be as in Theorem~\ref{auxilary result Ax^p+By^p+Cz^p=0 over K}, and let $(a,b,c)\in K^3$ be a non-trivial solution to equation~\eqref{Ax^p+By^p+Cz^p=0} with exponent $p>V$. Rescale $(a,b,c)$ so that $(a,b,c) \in \mcO_K^3$ and $G_{a,b,c}=\mfm$ for some $\mfm \in H$, and let $E'$ be the elliptic curve supplied by Theorem~\ref{auxilary result Ax^p+By^p+Cz^p=0 over K}. Since $E^\prime$ has good reduction away from $S_K^{\prime} \cup \{\mfm\}$ and potentially good reduction at $\mfm$, we get $j_{E^\prime} \in \mcO_{S_K^\prime}$.
	As in the proof of Theorem~\ref{main result1 for Ax^p+By^p+Cz^p=0}, choose a $\lambda$-invariant of $E'$ and put $\mu=1-\lambda$. Then $(\lambda,\mu)$ is a solution of the $S_K'$-unit equation~\eqref{S_K-unit solution}, and
	$
		j_{E'}=2^8\frac{(1-\lambda\mu)^3}{(\lambda\mu)^2}.
	$
	
	By the hypotheses of Theorem~\ref{main result2 for Ax^p+By^p+Cz^p=0}, there exists $\mfP\in U_K$ satisfying the two conditions in~\eqref{assumption for main result2 for Ax^p+By^p+Cz^p=0}, and the required coefficient conditions.
	Now, by using the same argument as in the proof of Theorem~\ref{main result1 for Ax^p+By^p+Cz^p=0}, 
	we conclude that $v_\mfP(j_{E^\prime})\geq0$ from the unit bound. Moreover,
	$
		v_\mfP(j_{E^\prime})\equiv8v_\mfP(2)-2v_\mfP(\lambda\mu)\pmod3.
	$
	Since $v_\mfP(\lambda\mu)\equiv v_\mfP(2)\pmod3$, we have
	$
		v_\mfP(j_{E^\prime})\equiv6v_\mfP(2)\equiv0\pmod3.
	$
	Therefore, $v_\mfP(j_{E^\prime}) \geq0$ and $3 \mid v_\mfP(j_{E^\prime})$, which contradicts Theorem~\ref{auxilary result Ax^p+By^p+Cz^p=0 over K}(4).
\end{proof}

Finally, we end this section with the following proposition.
\begin{prop}
	\label{relax the assumption}
	If $A,B,C \in \Z \setminus \{0\}$, then Theorem~\ref{main result2 for Ax^p+By^p+Cz^p=0} continues to hold over $K$ even without the assumption $A\pm B \pm C \neq 0$.
\end{prop}
\begin{proof}
	In order to prove Theorem~\ref{main result2 for Ax^p+By^p+Cz^p=0} without the condition $A\pm B \pm C \neq 0$, it suffices to prove the modularity result, namely Theorem~\ref{modularity of Frey curve x^2=By^p+Cz^p over W_K} without the condition $A\pm B \pm C \neq 0$.
	Using the same argument as in the proof of Theorem~\ref{modularity of Frey curve x^2=By^p+Cz^p over W_K} for $(a,b,c)\in K^3$, there are $\mu_1,\ldots,\mu_t\in K$ such that the Frey curve $E/K$ is modular whenever $\mu\notin\{\mu_1,\ldots,\mu_t\}$, where
	$
		j_E=2^8\frac{(\mu^2-\mu+1)^3}{\mu^2(\mu-1)^2},
	$ $
		\mu=-\frac{Bb^p}{Aa^p}.
	$
	We may assume that $\mu_1,\ldots,\mu_t\notin\Q^*$, because elliptic curves over $\Q$ are modular. If $\mu=\mu_k$, then
	$
		\left(\frac{b}{a}\right)^p=-\frac{A\mu_k}{B}.
	$
	Only finitely many prime exponents $p$ can satisfy this equation unless $-A\mu_k/B$ is a root of unity. In the exceptional case, since $K$ is totally real, $b/a=\pm1$ and hence $\mu_k=\pm B/A\in\Q^*$, a contradiction. Taking $D$ larger than the finitely many remaining exceptional exponents proves the required modularity statement.
\end{proof}

\section{Local criteria for solutions of \texorpdfstring{$Ax^p+By^p+Cz^p=0$}{Axp+Byp+Czp=0} over \texorpdfstring{$K$}{K}}
\label{section for loc criteria}
In this section, we give several purely local criteria under which equation~\eqref{Ax^p+By^p+Cz^p=0} has no asymptotic solution in $W_K$ or in $\mcO_K^3$. Throughout this section, we assume
$
	A,B,C\in\{2^ru:r\in\Z_{\geq0},\ u\in\mcO_K^*\};$ 
in particular, $S_K'=S_K$.
%
\subsection{Quadratic field}
\label{quadratic field}
First, we give criteria for quadratic fields under which Theorem~\ref{main result1 for Ax^p+By^p+Cz^p=0} applies. The following proposition is similar to \cite[Corollary 7.4]{KS24 Diop1}.

\begin{prop}
	\label{loc crit for real quadratic field over W_K}
	Let $d \geq 2$ be a square-free integer, let $K=\Q(\sqrt{d})$, and let $\mfP$ be the unique prime of $K$ above $2$. Assume that $\alpha_\mfP\beta_\mfP\gamma_\mfP\ne0$ and that one of the following conditions holds:
	\begin{enumerate}
		\item $d \equiv 3 \pmod 8$;
		\item $d \equiv 5 \pmod 8$;
		\item $d \equiv 6 \text{ or } 10 \pmod {16}$;
		\item $d \equiv 2 \pmod {16}$ and $d$ has some prime divisor $q \equiv 5 \text{ or } 7 \pmod 8$;
		\item $d \equiv 14 \pmod {16}$ and $d$ has some prime divisor $q \equiv 3 \text{ or } 5 \pmod 8$.
	\end{enumerate} 
	Then the generalized Fermat equation $Ax^p+By^p+Cz^p=0$ has no asymptotic solution in $W_K$.
\end{prop}

\begin{proof}
	Let $\mfP \in S_K$. By~\cite[Table 1 in \S6]{FS15}, every solution to the $S_K$-unit equation $\lambda + \mu=1$ satisfies $\max \left\{|v_\mfP(\lambda)|,|v_\mfP(\mu)| \right\} \leq 4 v_\mfP(2)$.
	Hence, the proof of the proposition follows from Theorem~\ref{main result1 for Ax^p+By^p+Cz^p=0}.
\end{proof}

We now give quadratic-field criteria under which Theorem~\ref{main result2 for Ax^p+By^p+Cz^p=0} applies. The following proposition is similar to \cite[Corollary 7.5]{KS24 Diop1}.
\begin{prop}
	\label{loc crit for real quadratic field over K}
	Let $d>6$ be square-free, and let $K=\Q(\sqrt d)$ satisfy one of the field conditions in Proposition~\ref{loc crit for real quadratic field over W_K}. Let $\mfP$ be the unique prime of $K$ above $2$. Suppose Conjecture~\ref{ES conj} holds over $K$. If $A\pm B\pm C\ne0$, $\max\{v_\mfP(A),v_\mfP(BC)\}\leq4v_\mfP(2)$, and $v_\mfP(ABC)\equiv0$ or $2v_\mfP(2)\pmod3$, then the equation $Ax^p+By^p+Cz^p=0$ has no asymptotic solution in $K^3$.
\end{prop}

\begin{proof}
	Since $[K: \Q]=2$, $S_K=U_K$. Again, by~\cite[Table 1 in \S6]{FS15}, the $S_K$-unit equation $\lambda + \mu=1$ has only the irrelevant solutions $(2,-1)$, $(-1,2)$, and $(\frac{1}{2}, \frac{1}{2})$. These satisfy $\max \left\{|v_\mfP(\lambda)|,|v_\mfP(\mu)| \right\}=v_\mfP(2)<4v_\mfP(2)$ and $v_\mfP(\lambda\mu)=v_\mfP(2)$ or $-2v_\mfP(2)$. Therefore, $v_\mfP(\lambda\mu)\equiv v_\mfP(2)\pmod3$.
	Hence, the proof of the proposition follows from Theorem~\ref{main result2 for Ax^p+By^p+Cz^p=0}.
\end{proof}
\subsection{Proof of Theorem~\ref{thm for density}}
\label{Proof of the thm for density}
In this subsection, we will prove Theorem~\ref{thm for density}. To prove this, we need to recall the absolute density of any subset $U \subseteq \N$.
\begin{dfn}
	For $S \subseteq \N$ and $x>0$, let $S(x):=\{d \in S : d \leq x\}$. The absolute density of $S$ is defined by
	$$ \delta(S):= \lim_{x \to \infty}\frac{\# S(x)}{x},$$
	if the limit exists.
\end{dfn}
The following theorem is useful in the proof of Theorem~\ref{thm for density}.
\begin{thm}\rm(\cite[Theorem 10]{FS15})
	\label{FS thm for density}
	For $r\in\Z$ and $N\in\N$, let $\N_{r,N}^{\text{sf}}:=\{d\in\N^{\text{sf}}:d\equiv r\pmod N\}$. If $s:=\gcd(r,N)$ is square-free, then
	\[
		\#\N_{r,N}^{\text{sf}}(x)
		\sim
		\frac{\varphi(N)}{s\varphi(N/s)N\prod_{q\mid N}(1-1/q^2)}
		\cdot\frac{6}{\pi^2}x,
	\]
	where $\varphi$ denotes Euler's totient function.
\end{thm}
The following lemma is very useful for the proof of Theorem~\ref{thm for density}.
\begin{lem}\rm(\cite[Lemma 7.1]{FS15})
	\label{FS lem for density}
	Let $V:=\{d \in \N^{\text{sf}} :  \text{ the equation } \lambda+\mu=1, \ \lambda, \mu \in \mcO_{S_K}^\ast  \text{has only irrelevant solutions for } K=\Q(\sqrt{d}) \}$. Then $\delta(\N^{\text{sf}} \setminus V)=0$.
\end{lem}
We are now ready to prove Theorem~\ref{thm for density}
\begin{proof}[Proof of Theorem~\ref{thm for density}]
	Note that $\N^{\text{sf}}=\N_{0,1}^{\text{sf}}$. By Theorem~\ref{FS thm for density},
	$
		\#\N^{\text{sf}}(x)=\#\N_{0,1}^{\text{sf}}(x)\sim\frac{6}{\pi^2}x,
	$
	hence $\delta(\N^{\text{sf}})=6/\pi^2$. For any $S\subseteq\N^{\text{sf}}$, the density $\delta(S)$ exists if and only if $\delta_{\text{rel}}(S)$ exists; in that case,
	$
		\delta_{\text{rel}}(S)=\frac{\delta(S)}{\delta(\N^{\text{sf}})}
		=\frac{\pi^2}{6}\delta(S).
	$
	By definition, $U=V$. Lemma~\ref{FS lem for density} gives $\delta(\N^{\text{sf}}\setminus U)=0$, and therefore $\delta_{\text{rel}}(\N^{\text{sf}}\setminus U)=0$. Hence $\delta_{\text{rel}}(U)=1$. For $d\in U$, the three irrelevant solutions satisfy
	$
		\max\{|v_\mfP(\lambda)|,|v_\mfP(\mu)|\}=v_\mfP(2),
 $ and $
		v_\mfP(\lambda\mu)\equiv v_\mfP(2)\pmod3
	$
	for every $\mfP\in S_K$. The final assertions now follow from Theorems~\ref{main result1 for Ax^p+By^p+Cz^p=0} and~\ref{main result2 for Ax^p+By^p+Cz^p=0}.
\end{proof}

\subsection{Odd degree}
First, we give odd-degree criteria under which Theorem~\ref{main result1 for Ax^p+By^p+Cz^p=0} applies. The following proposition is similar to \cite[Corollary 7.3]{KS24 Diop1}.
\begin{prop}
	\label{loc crit for odd degree field over W_K}
	Let $n=[K:\Q]$ be odd, and assume that $2$ is either inert or totally ramified in $K$. Let $\mfP$ be the unique prime of $K$ above $2$, and assume $\alpha_\mfP\beta_\mfP\gamma_\mfP\ne0$. Suppose one of the following conditions holds:
	\begin{enumerate}
		\item there is a prime $l>5$ such that $(n,l-1)=1$ and $l$ is totally ramified in $K$;
		\item $3$ splits completely in $K$.
	\end{enumerate}
	Then the equation $Ax^p+By^p+Cz^p=0$ has no asymptotic solution in $W_K$.
\end{prop}

\begin{proof}
	Let $\mfP \in S_K$ be the unique prime ideal lying above $2$.
	Now, arguing as in the proof of \cite[Corollary 7.3]{KS24 Diop1}, we find that every solution $(\lambda, \mu)$ of the $S_K$-unit equation $\lambda + \mu=1$ satisfies 	$\max \left\{|v_\mfP(\lambda)|,|v_\mfP(\mu)| \right\}<2 v_\mfP(2)$. Hence, the proof of the proposition follows from Theorem~\ref{main result1 for Ax^p+By^p+Cz^p=0}.
\end{proof}

We now give odd-degree criteria under which Theorem~\ref{main result2 for Ax^p+By^p+Cz^p=0} applies. The following proposition is similar to \cite[Propositions 7.1, 7.2]{KS24 Diop1}.

\begin{prop}
	\label{local criteria1 for odd defree over K}
	\label{local criteria2 for odd defree over K}
	Let $n=[K:\Q]$ be odd. Suppose one of the following conditions holds:
	\begin{enumerate}
		\item there is a prime $l>5$ such that $(n,l-1)=1$, $l$ is totally ramified in $K$, and $2$ is inert in $K$;
		\item $3\nmid n$, $2$ is inert in $K$, and $3$ splits completely in $K$.
	\end{enumerate}	
	Let $\mfP$ be the unique prime of $K$ above $2$. Assume that $A\pm B\pm C\ne0$ and
	\[
		\max\{v_\mfP(A),v_\mfP(BC)\}\leq4v_\mfP(2),
		\qquad
		v_\mfP(ABC)\equiv0\ \text{or}\ 2v_\mfP(2)\pmod3.
	\]
	Then the equation $Ax^p+By^p+Cz^p=0$ has no asymptotic solution in $K^3$.
\end{prop}

\begin{proof}
	Arguing as in the proof of \cite[Propositions 7.1, 7.2]{KS24 Diop1}, we find that every solution $(\lambda, \mu)$ of the $S_K$-unit equation $\lambda + \mu=1$ satisfies $\max \left\{|v_\mfP(\lambda)|,|v_\mfP(\mu)| \right\}= v_\mfP(2)$.
	Let $s:=\max \left\{|v_\mfP(\lambda)|,|v_\mfP(\mu)| \right\}$. Then $s =v_\mfP(2)>0$. Since $\lambda +\mu =1$, we get $v_\mfP(\lambda\mu)= -2s$ or $s$. Therefore, $v_\mfP(\lambda\mu) \equiv s= v_\mfP(2) \pmod 3$. 
	Hence, the proof of the proposition follows from Theorem~\ref{main result2 for Ax^p+By^p+Cz^p=0}.
\end{proof}

\section*{Acknowledgments} 
The author would like to express his sincere gratitude to Prof. Nuno Freitas for his invaluable assistance in understanding the article \cite{FS15}. The author also extends his heartfelt appreciation to Prof. Narasimha Kumar for the insightful discussions and helpful comments. The author is grateful to the anonymous referee for the thoughtful comments and valuable suggestions, which have greatly improved the quality of this article.

\end{document}